\documentclass[12pt]{amsart}
\newtheorem{theorem}{Theorem}
\newtheorem{corollary}[theorem]{Corollary}
\newenvironment{mat}{\footnotesize$\left\lgroup\begin{array}{ccc}}
                     {\end{array}\right\rgroup$}
\newcommand{\spn}{\operatorname{span}}
\newcommand{\md}{\operatorname{mod}}
\newcommand{\ddt}{\mbox{\small$\displaystyle\frac{\partial}{\partial t}$}}
\newcommand{\ddx}{\mbox{\small$\displaystyle\frac{\partial}{\partial x}$}}
\newcommand{\ddy}{\mbox{\small$\displaystyle\frac{\partial}{\partial y}$}}
\begin{document}
\title{Preferred Parameterisations on Homogeneous Curves}
\author[Michael Eastwood]{Michael Eastwood}
\address{Department of Pure Mathematics\\ University of Adelaide,\newline
         \indent South AUSTRALIA 5005}
\email{meastwoo@maths.adelaide.edu.au}
\author[Jan Slov\'ak]{Jan Slov\'ak}
\address{Department of Algebra and Geometry\\ Masaryk University,\newline
\indent Jan\'a\v ckovo N\'am 2a\\ 662 95 Brno\\ CZECH REPUBLIC}
\email{slovak@math.muni.cz}
\thanks{Support from the Australian Research Council is gratefully
acknowledged. The second author was partially supported by GACR 201/02/1390}
\begin{abstract} We show how to specify preferred parameterisations on
a homogeneous curve in an arbitrary homogeneous space. We apply these results
to limit the natural parameters on distinguished curves in parabolic
geometries.
\end{abstract}
\maketitle
\section{Introduction}
This article is motivated by the theory of {\em distinguished curves} in
{\em parabolic geometries}, as developed in~\cite{csz}. A parabolic geometry
is, by definition, modelled on a homogeneous space of the form $G/P$ where $G$
is a real semisimple Lie group and $P$ is a parabolic subgroup. (There is also
a complex theory
which corresponds to the choices of complex $G$'s and $P$'s with specific
curvature restrictions for the holomorphic cases.)
The notion of Cartan connection
replaces the Maurer-Cartan form on~$G$, viewed as a principal fibre bundle over
$G/P$ with structure group~$P$, and much of the geometry of $G/P$ automatically
carries over to parabolic geometries in general (see also~\cite{sharpe}). In
particular, the curves on $G/P$ obtained by exponentiating elements in the Lie
algebra ${\mathfrak g}$ of $G$ have counterparts in general obtained by
{\em development} under the Cartan connection. These matters are thoroughly
discussed in \cite{csz} and will not be repeated here. Suffice it to say that
results concerning distinguished curves on $G/P$ have immediate consequences
for the corresponding general parabolic geometry. Here, we shall discuss only
the homogeneous setting~$G/P$.

\section{Generalities on $G/P$}\label{generalities}
Firstly, let us discuss a general homogeneous space, namely a smooth manifold
$M$ equipped with the smooth transitive action of a real Lie group~$G$. Each
$X\in{\mathfrak g}$ gives a 1-parameter Lie subgroup $t\mapsto\exp(tX)$ of $G$
and hence to a parameterised curve $t\mapsto\exp(tX)m$ through $m\in M$, which
we shall suppose to be non-constant. Conversely, without the parameterisation,
such a curve is {\em homogeneous}, namely it is the orbit of a Lie subgroup of
the symmetry group~$G$.

To investigate homogeneous curves on $M$ we may as well choose a basepoint
$m_\circ\in M$ and consider only curves passing through~$m_\circ$. All other
homogeneous curves are obtained by translation under the action of~$G$. Let $P$
denote the stabiliser subgroup of~$m_\circ$ so that $M=G/P$. We shall now
suppose that $G$ is semisimple and $P$ is parabolic. In this case, there is a
splitting
$${\mathfrak g}={\mathfrak p}\oplus{\mathfrak n}$$
into subalgebras with ${\mathfrak n}$ nilpotent (as in~\cite{csz}). This
splitting is not canonical. It is, however, well-defined up to the Adjoint
action of $P$ and we obtain, therefore, a preferred subset
\begin{equation}\label{preferredsubset}
\{{\mbox{Ad}}_pX\mbox{ s.t. }p\in P\mbox{ and }X\in{\mathfrak n}\}
\subset{\mathfrak g},\end{equation}
which we may use to generate homogeneous curves. Such curves are evidently
non-constant but not all non-constant homogeneous curves arise in this way.
These special curves are said to be {\em distinguished}. Equivalently,
distinguished curves through the basepoint $m_\circ\in M$ are those of the form
$t\mapsto p\exp(tX)m_\circ$ for some $p\in P$ and $X\in{\mathfrak n}$.

As an example, consider $G={\mathrm{SL}}(3,{\mathbb R})$ with $P$ the upper
triangular matrices. We may take
\begin{equation}\label{lowertriangular}{\mathfrak n}=
\left\{\mbox{\begin{mat}0&0&0\\ *&0&0\\ *&*&0\end{mat}}\right\}.\end{equation}
Then
$$t\mapsto\mbox{\begin{mat}1&1&0\\ 0&1&0\\ 0&0&1\end{mat}}
\exp\left(t\mbox{\begin{mat}0&0&0\\ 1&0&0\\ 0&0&0\end{mat}}\right)\md P
=\mbox{\begin{mat}1+t&1&0\\ t&1&0\\ 0&0&0\end{mat}}\md P$$
and
$$t\mapsto
\exp\left(t\mbox{\begin{mat}0&0&0\\ 1&0&0\\ 1&1&0\end{mat}}\right)\md P
=\mbox{\begin{mat}1&0&0\\ t&1&0\\ t+\frac12t^2&t&1\end{mat}}\md P$$
are typical distinguished curves whereas
$$t\mapsto
\exp\left(t\mbox{\begin{mat}0&-1&0\\ 1&0&0\\ 0&0&0\end{mat}}\right)\md P
=\mbox{\begin{mat}\cos t&-\sin t&0\\ \sin t&\cos t&0\\ 0&0&1\end{mat}}\md P$$
is homogeneous but (with this parameterisation) not distinguished.

Suppose $y\mapsto\gamma(t)\in M$ is a distinguished curve with
$\gamma(0)=m_\circ$ and let $C$ denote its unparameterised image. In this
article, we shall answer the question `what are the possible
reparameterisations of $C$ as a distinguished curve?' A direct approach to this
question is given in~\cite[\S3]{csz}. Here, we shall reason indirectly by
firstly establishing the following on general grounds.

\begin{theorem}\label{maintheorem} Let $C$ be an unparameterised distinguished
curve passing through $m_\circ\in M=G/P$. The freedom in reparameterising $C$
with origin at $m_\circ$ is of two possible types:--
\begin{flushleft}\quad\begin{tabular}{lll}
{\bf affine}&$t\mapsto at$&for $a\not=0$\\
{\bf projective}&$t\mapsto at/(ct+1)$&for $a\not=0$ and $c$ arbitrary.
\end{tabular}\end{flushleft}
\end{theorem}
\noindent If we drop the requirement that the parameter be zero at~$m_\circ$,
then translation is also allowed so the freedom becomes
$$t\mapsto at+b\quad\mbox{or}\quad t\mapsto\frac{at+b}{ct+d},$$
respectively. The proof of Theorem~\ref{maintheorem} will be given
in~\S\ref{mainsection}. Once this theorem is established, it is a matter of
elementary computation to decide, for a given~$C$, which type of freedom
pertains. Examples will be given in~\S\ref{mainsection}. For the proof of
Theorem~\ref{maintheorem} we shall need some general considerations as in the
following section.

\section{Lie algebras of vector fields in one dimension}
The following is a classical topic and Theorem~\ref{local} is due to
Lie~\cite{lie} (see also~\cite{strigunova}). We are grateful to Ian Anderson
for pointing out to us the translation and commentary on Lie's article given by
Ackerman and Hermann~\cite{ah}. Nevertheless, we believe that it is useful to
given an independent, elementary, and self-contained treatment.

\begin{theorem}\label{one} Suppose ${\mathfrak g}$ is a finite-dimensional
subalgebra of the Lie algebra of smooth vector fields on~${\mathbb R}$. Let $x$
be the standard co\"ordinate on ${\mathbb R}$ and suppose
${\mathfrak g}\ni\partial/\partial x$. Then ${\mathfrak g}$ is one of the
following:--
$$\begin{array}c\displaystyle
{\mathfrak g}=\spn\left\{\ddx\right\}\quad
{\mathfrak g}=\spn\left\{\ddx,e^{\lambda x}\ddx\right\}\quad
{\mathfrak g}=\spn\left\{\ddx,x\ddx\right\}\\[12pt] \displaystyle
{\mathfrak g}=\spn\left\{\ddx,\sin(\lambda x)\ddx,
              \cos(\lambda x)\ddx\right\}\quad
{\mathfrak g}=\spn\left\{\ddx,x\ddx,x^2\ddx\right\}\\[12pt] \displaystyle
{\mathfrak g}=\spn\left\{\ddx,\sinh(\lambda x)\ddx,
              \cosh(\lambda x)\ddx\right\}.
\end{array}$$
\end{theorem}
\begin{proof}If $\dim{\mathfrak g}=1$, then
${\mathfrak g}=\spn\{\partial/\partial x\}$ are we are done. Next, if
$\dim{\mathfrak g}=2$, then
${\mathfrak g}=\spn\{\partial/\partial x,g(x)\partial/\partial x\}$ for some
smooth non-constant function~$g(x)$. Now,
$$\left[\ddx,g(x)\ddx\right]=g^\prime(x)\ddx$$
so closure under Lie bracket implies $g^\prime(x)=\mu+\lambda g(x)$. This is a
differential equation we may solve:--
$$\begin{array}{rcll}
g(x)&=&Ce^{\lambda x}+D&\mbox{if }\lambda\not=0\\[3pt]
\mbox{or}\quad g(x)&=&\mu x+C&\mbox{if }\lambda=0.\end{array}$$
These are the two-dimensional subalgebras stated in the theorem.

Now suppose $\dim{\mathfrak g}=k+1\geq 3$ and choose a basis
$$\ddx,g_1(x)\ddx,\ldots,g_k(x)\ddx$$
of ${\mathfrak g}$. {From} closure under Lie bracket of $\partial/\partial x$
with the other basis vectors, we immediately encounter a system of ordinary
differential equations with constant coefficients
$$g_i^\prime(x)=\mu_i+\sum_{j=1}^k\lambda_{ij}g_j(x),
\quad\mbox{for }i=1,\ldots,k.$$
We may conclude that the functions $g_i(x)$ and, therefore, all vector fields
in ${\mathfrak g}$ are real-analytic.

Since $\dim{\mathfrak g}\geq 3$, there is a vector field
$g(x)\partial/\partial x\in{\mathfrak g}$ with
$$g(x)=x^N+\cdots\quad\mbox{for some }N\geq 2.$$
Because ${\mathfrak g}$ is finite-dimensional, we may choose $g(x)$ with $N$
maximal. But then
$$\begin{array}{rcl}\displaystyle
{\mathfrak g}\ni\left[\left[\ddx,g(x)\ddx\right],g(x)\ddx\right]&=&
\left[g^\prime(x)\ddx,g(x)\ddx\right]\\[12pt]
&=&\big((g^\prime(x))^2-g(x)g^{\prime\prime}(x)\big)\ddx\\[8pt]
&=&\big(Nx^{2(N-1)}+\cdots\big)\ddx,
\end{array}$$
contradicting maximality of $N$ unless $N=2$. Therefore,
$\dim{\mathfrak g}=3$ and
\begin{equation}\label{thisistheLiealgebra}
{\mathfrak g}=\spn\left\{\ddx,g(x)\ddx,g^\prime(x)\ddx\right\},
\end{equation}
where
\begin{equation}\label{thisisthefunction}
g(x)= x^2+ax^3+\cdots.\end{equation}
But then ${\mathfrak g}$ contains the vector field
$$\begin{array}{rcl}\left[g^\prime(x)\ddx,g(x)\ddx\right]-2g(x)\ddx
&\!=\!&\big((g^\prime(x))^2-g(x)g^{\prime\prime}(x)-2g(x)\big)\ddx\\[8pt]
&\!=\!&\big(2ax^3+\cdots\big)\ddx,
\end{array}$$
again contradicting maximality of $N$ unless $a=0$. Now, in order for
(\ref{thisistheLiealgebra}) to be closed under Lie bracket we must have
$$g^{\prime\prime}(x)=\left[\ddx,g^\prime(x)\ddx\right]
\in\spn\left\{\ddx,g(x)\ddx,g^\prime(x)\ddx\right\}$$
and to be, in addition, consistent with $a=0$ in (\ref{thisisthefunction}), we
conclude that
$$g^{\prime\prime}(x)=2+\nu g(x),\quad\mbox{for some constant }\nu.$$
This differential equation, with initial conditions imposed by
(\ref{thisisthefunction}), has solutions
$$\begin{array}{rcll}
g(x)&=&(2/\lambda^2)(\cos(\lambda x)-1)&\mbox{if }\nu<0\\[3pt]
\mbox{or}\quad g(x)&=&x^2&\mbox{if }\nu=0\\[3pt]
\mbox{or}\quad g(x)&=&(2/\lambda^2)(\cosh(\lambda x)-1)&\mbox{if }\nu>0.
\end{array}$$
It remains to observe that (\ref{thisistheLiealgebra}) is, indeed, closed under
Lie bracket in these cases.
\end{proof}
Notice that this proof is local: the same conclusion holds for vector fields on
any open interval $(a,b)\subseteq{\mathbb R}$. Globally on ${\mathbb R}$, the
various subalgebras given in the statement of Theorem~\ref{one} are distinct.
Locally, however, this distinction evaporates leaving only the dimension. There
is a co\"ordinate change near the origin:--
$$y=\frac{1-e^{-\lambda x}}{\lambda}\quad\Rightarrow\quad
e^{\lambda x}\ddx=\ddy\quad\mbox{and}\quad\ddx=(1-\lambda y)\ddy$$
whence
$$\spn\left\{\ddx,e^{\lambda x}\ddx\right\}\cong
\spn\left\{y\ddy,\ddy\right\}\cong\spn\left\{\ddx,x\ddx\right\}$$
whilst the co\"ordinate change $y=\tan((\lambda x)/2)$ gives
$$\ddx=\frac{\lambda}{2}(1+y^2)\ddy,\;
\sin(\lambda x)\ddx=\lambda y\ddy,\;
\cos(\lambda x)\ddx=\frac{\lambda}{2}(1-y^2)\ddy$$
whence
$$\spn\left\{\ddx,\sin(\lambda x)\ddx,\cos(\lambda x)\ddx\right\}\cong
\spn\left\{\ddx,x\ddx,x^2\ddx\right\}$$
and $y=\tanh((\lambda x)/2)$ gives
$$\ddx=\frac{\lambda}{2}(1-y^2)\ddy,\;
\sinh(\lambda x)\ddx=\lambda y\ddy,\;
\cosh(\lambda x)\ddx=\frac{\lambda}{2}(1+y^2)\ddy$$
whence
$$\spn\left\{\ddx,\sinh(\lambda x)\ddx,\cosh(\lambda x)\ddx\right\}\cong
\spn\left\{\ddx,x\ddx,x^2\ddx\right\}.$$
We have proved the following.
\begin{theorem}\label{local}
Suppose ${\mathfrak s}$ is a finite-dimensional subalgebra of
the Lie algebra of vector fields in a neighbourhood of the origin
in~${\mathbb{R}}$. Suppose ${\mathfrak s}$ contains a vector field that does
not vanish at the origin. Then there is a neighbourhood $\,U$ of the origin and
a change of co\"ordinates such that one of the following three possibilities
holds.
\begin{equation}\label{listed}\begin{array}c\displaystyle
{\mathfrak s}|_U\cong\spn\left\{\ddx\right\}\qquad
{\mathfrak s}|_U\cong\spn\left\{\ddx,x\ddx\right\}\\[12pt]\displaystyle
{\mathfrak s}|_U\cong\spn\left\{\ddx,x\ddx,x^2\ddx\right\}.
\end{array}\end{equation}
\end{theorem}

\section{Reparameterisations}\label{mainsection}
Let $C$ be an arbitrary smooth connected curve in a smooth manifold $M$
homogeneous under the action $\rho:G\times M\to M$ of a connected Lie
group~$G$. There is a homomorphism of Lie algebras $\dot\rho:{\mathfrak
g}\to\mbox{Vect}(M)$ given by
$$\dot\rho(X)(m)=\ddt\big(\exp(-tX)m\big)|_{t=0}$$
and the {\em symmetry algebra} of $C$ is defined by
$${\mathfrak s}=\{X\in{\mathfrak g}\mbox{ s.t. }\dot\rho(X)(m)\mbox{ is
tangent to }C\mbox{ for all }m\in C\}.$$
Clearly, ${\mathfrak s}$ is a subalgebra of ${\mathfrak g}$ and $C$ is
homogeneous if and only if $\dot\rho({\mathfrak s})|_C$ contains non-trivial
vector fields at each point of~$C$. In this
case, we may invoke Theorem~\ref{local} to conclude that
$\dot\rho({\mathfrak s})|_C$ is at most three-dimensional and locally has one
of the three forms listed in~(\ref{listed}).

Now suppose that $C$ is homogeneous and pick a basepoint $m_\circ\in C$.
Suppose that $X\in{\mathfrak s}\subset{\mathfrak g}$ is nilpotent in
${\mathfrak g}$ and $\dot\rho(X)(m_\circ)\not=0$. Then we shall say that
$$t\mapsto\exp(tX)m_\circ\in C$$
is a {\em preferred parameterisation} of~$C$.
\begin{theorem}\label{freedom} The freedom in reparameterising a homogeneous
curve with a preferred parameter is one of two possible types:--
\begin{flushleft}\quad\begin{tabular}{lll}
{\bf affine}&$t\mapsto at$&for $a\not=0$\\
{\bf projective}&$t\mapsto at/(ct+1)$&for $a\not=0$ and $c$ arbitrary.
\end{tabular}\end{flushleft}
\end{theorem}
\begin{proof} Since $X$ is nilpotent in ${\mathfrak g}$, certainly
$\dot\rho(X)$ is nilpotent in~$\dot\rho({\mathfrak s})|_C$. By inspection, we
may find the nilpotent elements in each of the local forms~(\ref{local}):--
$$\begin{array}{rcl}
\displaystyle\spn\left\{\ddx\right\}&\ni&a\ddx\\[12pt]
\displaystyle\spn\left\{\ddx,x\ddx\right\}&\ni&a\ddx\\[12pt]
\displaystyle\spn\left\{\ddx,x\ddx,x^2\ddx\right\}&\ni&(p-qx)^2\ddx.
\end{array}$$
In the first two cases,
$$a\ddx=\ddt\iff x=at,$$
which gives affine freedom, whilst in the third case
$$(p-qx)^2\ddx=\ddt\iff x=\frac{p^2t}{1+pqt},$$
which gives projective freedom.
\end{proof}
\noindent{\em Proof of Theorem~\ref{maintheorem}.} The
parameterisations on a distinguished curve have the form
$t\mapsto\exp(tY)m_\circ$ where $Y$ is $P$-conjugate to an element of
${\mathfrak n}$ in accordance with~(\ref{preferredsubset}). Certainly, there is
affine freedom in such a parameterisation because $Y$ can be replaced by~$aY$.
But the allowed $Y$ are, in particular, nilpotent. Therefore, the
parameterisations on $C$ as a distinguished curve are {\em ipso facto}
preferred parameterisations on $C$ as a homogeneous curve.
Theorem~\ref{freedom} now implies that, if there is any additional freedom, it
must be projective. But just one projective transformation, together with
affine freedom, generates all projective freedom and the proof is complete.
\hfill$\qed$

Theorem~\ref{maintheorem} is useful in practice. Consider the general
distinguished curve $t\mapsto p\exp(tX)m_\circ$ for fixed $p\in P$ and
$X\in{\mathfrak n}$. The dichotomy offered by Theorem~\ref{maintheorem} implies
that if there are reparameterisations other than affine, then the specific
projective freedom $t\mapsto t/(t+1)$ occurs. In this case, we may find
$q\in P$ and $Y\in{\mathfrak n}$ such that
$$p\exp(tX)=q\exp\big(\frac{t}{t+1}Y\big)\md P,\quad\forall t$$
or, equivalently,
\begin{equation}\label{criterion}
\exp\big(-\frac{t}{t+1}Y\big)r\exp(tX)\in P,\quad\forall t\end{equation}
where $r=q^{-1}p\in P$. The existence of suitable $r\in P$ and
$Y\in{\mathfrak n}$ is a restriction on~$X$. Furthermore, if we adopt the Levi
decomposition $P=LU$ corresponding to our choice of ${\mathfrak n}$, then the
$L$-component of $r$ may be absorbed into~$Y$. Hence, Theorem~\ref{maintheorem}
has the following corollary.
\begin{corollary}\label{maincorollary}
Suppose $P=LU$ is a Levi decomposition of a parabolic subgroup $P$ of a
semisimple Lie group~$G$. Let ${\mathfrak g}={\mathfrak p}\oplus{\mathfrak n}$
be the associated decomposition of the Lie algebra of~$G$. Then the
distinguished curve $t\mapsto p\exp(tX)\md P$ admits a projective
reparameterisation if and only if there are $r\in U$ and $Y\in{\mathfrak n}$
such that {\rm{(\ref{criterion})}} holds.
\end{corollary}

We close this article with a complete analysis of the distinguished curves in
the real flag manifold ${\mathrm{SL}}(3,{\mathbb R})/P$ where $P$ is the
subgroup of upper triangular matrices. As already remarked in
\S\ref{generalities}, we may take ${\mathfrak n}$ to be the strictly lower
triangular matrices~(\ref{lowertriangular}). We shall use
Corollary~\ref{maincorollary} with $U$ taken to be the upper triangular
matrices with $1$'s along the diagonal. Consider, for example, the
distinguished curve
\begin{equation}\label{curve}
t\mapsto\mbox{\begin{mat}1&1&0\\ 0&1&0\\ 0&0&1\end{mat}}
\exp\left(t\mbox{\begin{mat}0&0&0\\ 1&0&0\\ 0&0&0\end{mat}}\right)\md P.
\end{equation}
According to Corollary~\ref{maincorollary}, it admits a projective
reparameterisation if and only if we can find $a,b,c,u,v,w$ such that
$$\exp\left(-\frac{t}{t+1}
\mbox{\begin{mat}0&0&0\\ u&0&0\\ v&w&0\end{mat}}\right)
\mbox{\begin{mat}1&a&b\\ 0&1&c\\ 0&0&1\end{mat}}
\mbox{\begin{mat}1&0&0\\ t&1&0\\ 0&0&1\end{mat}}=
\mbox{\begin{mat}*&*&*\\ 0&*&*\\ 0&0&*\end{mat}}.$$
Multiplying through by $(t+1)^2$ yields
$$\mbox{\begin{mat}(t+1)^2&0&0\\ -t(t+1)u&(t+1)^2&0\\
        -t(t+1)v+\frac12t^2uw&-t(t+1)w&(t+1)^2\end{mat}}
\mbox{\begin{mat}1&a&b\\ 0&1&c\\ 0&0&1\end{mat}}
\mbox{\begin{mat}1&0&0\\ t&1&0\\ 0&0&1\end{mat}}$$
for the left hand side. Expanding and equating coefficients of $t$ to zero in
the subdiagonal entries gives algebraic equations for $a,b,c,u,v,w$ whose
general solutions are
$$\mbox{\begin{mat}0&0&0\\ u&0&0\\ v&w&0\end{mat}}=
\mbox{\begin{mat}0&0&0\\ 1&0&0\\ 0&0&0\end{mat}}\qquad
\mbox{\begin{mat}1&a&b\\ 0&1&c\\ 0&0&1\end{mat}}=
\mbox{\begin{mat}1&1&b\\ 0&1&c\\ 0&0&1\end{mat}}.$$
The existence of solutions shows that the distinguished curve (\ref{curve})
admits projective reparameterisations. On the other hand, this same exercise
for the curve
$$t\mapsto
\exp\left(t\mbox{\begin{mat}0&0&0\\ 1&0&0\\ 1&1&0\end{mat}}\right)\md P$$
gives an inconsistent set of equations for $a,b,c,u,v,w$. According to
Theorem~\ref{maintheorem} and Corollary~\ref{maincorollary}, it admits only
affine reparameterisations.

Notice that the criterion (\ref{criterion}) of Corollary~\ref{maincorollary}
depends only on the $L$-conjugacy class of $X\in{\mathfrak n}$. Therefore, to
say which distinguished curves admit projective reparameterisations it suffices
to say whether (\ref{criterion}) is satisfied for $X\in{\mathfrak n}$
normalised under the Adjoint action of~$L$.
We obtain the following table of
normal forms.
\smallskip\begin{center}\begin{tabular}{|c|c|}\hline
Normal form & Reparameterisation\\ \hline\hline
\mbox{\begin{mat}0&0&0\\ 1&0&0\\ 0&0&0\end{mat}} & projective\\ \hline
\mbox{\begin{mat}0&0&0\\ 0&0&0\\ 0&1&0\end{mat}} & projective\\ \hline
\mbox{\begin{mat}0&0&0\\ 1&0&0\\ 0&1&0\end{mat}} & projective\\ \hline
\mbox{\begin{mat}0&0&0\\ 1&0&0\\ 1&0&0\end{mat}} & projective\\ \hline
\mbox{\begin{mat}0&0&0\\ 0&0&0\\ 1&1&0\end{mat}} & projective\\ \hline
\mbox{\begin{mat}0&0&0\\ 0&0&0\\ 1&0&0\end{mat}} & projective\\ \hline
\mbox{\begin{mat}0&0&0\\ 1&0&0\\ x&1&0\end{mat}} & affine if $x\not=0$\\ \hline
\end{tabular}\end{center}\smallskip
We have to be careful, however, with the decision which of the above normal
forms give rise to different distinguished curves. In our case, the lines
four through six in the table are in the same orbit of the Adjoint action of
the entire $P$ and so the distinguished curves indicated by these lines
coincide. Indeed, a simple check reveals
$$
(\exp Z)^{-1}
\mbox{\begin{mat}0&0&0\\ 0&0&0\\ 1&0&0\end{mat}}
\exp Z
=\mbox{\begin{mat}0&0&0\\ 1&0&0\\ 1&0&0\end{mat}},\quad
Z= \mbox{\begin{mat}0&0&0\\ 0&0&-1\\ 0&0&0\end{mat}}
$$
while the other case is symmetric.
We should also like to remark, that the latter
observation yields a sufficient condition for coincidences of classes of
distinguished curves. There are also examples of such a coincidence where
the corresponding $L$--orbits are not in the same orbit of $P$. In our case,
however, the first three lines and the last two lines in the table obviously
produce different curves.

This completes our analysis of distinguished curves in this real flag
manifold. It is more efficient than the direct approach of~\cite{csz} because
Theorem~\ref{maintheorem} tells us, in advance, what sort of reparameterisation
we may expect on a distinguished curve.

\end{document}